\documentclass[12pt]{amsart}
\usepackage{amsmath,amsfonts,amsthm,amscd,lastpage,times}
\newtheorem{theorem}{Theorem}

\newtheorem{proposition}{Proposition}
\newtheorem{lemma}{Lemma}
\newtheorem{definition}{Definition}
\newtheorem{corollary}{Corollary}
\newtheorem{remark}{Remark}
\newtheorem{example}{Example}
\numberwithin{equation}{section}
\numberwithin{lemma}{section}
\numberwithin{claim}{section}
\input{diagrams.tex}
\newcommand{\bull}{\ensuremath{{}\bullet{}}}

\newcommand{\lam}{\ensuremath{\lambda}}

\newcommand{\cpn}{\ensuremath{\mathbb{P}^{N}}}
\newcommand{\slnc}{\ensuremath{SL(N+1,\mathbb{C})}}

\newcommand{\dlb}{\ensuremath{\overline{\partial}}}
\newcommand{\dl}{\ensuremath{\partial}}
\newcommand{\ra}{\ensuremath{\longrightarrow}}
\newcommand{\ba}{\ensuremath{\begin{align*}}}
\newcommand{\ea}{\ensuremath{\end{align*}}}

\newcommand{\vplt}{\ensuremath{\varphi_{\lambda(t)}}}
\newcommand{\vp}{\ensuremath{\varphi}}

\begin{document}
\title[CM Stability II]{CM Stability and the Generalized Futaki Invariant II}
\author{Sean Timothy Paul}
 \email{stpaul@math.wisc.edu}
 \address{Mathematics Department at the University of Wisconsin, Madison \newline 480 Lincoln Drive Madison Wi 53706 }
 \thanks{Supported by an NSF DMS grant \# 0505059}
\author{Gang Tian}
\email{tian@math.princeton.edu}
\address{Department of Mathematics, Fine Hall, Princeton University \newline Washington Road Princeton N.J. 08544}
\thanks{Supported by an NSF DMS grant}
\subjclass[2000]{Primary: 53C55; Secondary: 14D06}
 \keywords{Futaki Invariant, Stability, Mabuchi Energy, CM polarization . }

\date{March 31, 2008}
\maketitle
\vspace{-5mm}

 {\begin{abstract}{The Mabuchi K-energy map is exhibited as a singular metric on the refined CM polarization of any equivariant family $\mathbf{X}\overset{p}{\rightarrow} S$. Consequently we show that the generalized Futaki invariant is the leading term in the asymptotics of the reduced K-energy of the generic fiber of the map $p$. Properness of the K-energy implies that the generalized Futaki invariant is strictly negative.}
\end{abstract}
\section{Introduction}
\subsection{Statement of results}
Throughout this paper $\mathbf{X}$ and $S$ denote smooth, proper complex projective varieties  satisfying the following conditions.
 \smallskip
  \begin{enumerate}
  \item $\mathbf{X}\subset S\times\cpn $; \  $\cpn$  denotes the complex projective space of \emph{lines} in $\mathbb{C}^{N+1}$ .\\
  \item $p:=p_1: \mathbf{X}\rightarrow S$ is flat of relative dimension $n$, degree $d$ with Hilbert polynomial $P$.\\
  \item $L|_{\mathbf{X}_z}$ is very ample and the embedding ${\mathbf{X}_z}:=p_1^{-1}(z)\overset{L}{\hookrightarrow} \cpn$ is given by a complete linear system for $z\in S$.\\
\item There is an action of $G:=\slnc$ on the data compatible with the projection and the standard action on $\cpn$.
  \end{enumerate}

 \smallskip

  It is well known that  $(1)$  and (3) imply that 
  \begin{align}
  \mathbb{P}({p_1}_*L)\cong S\times\cpn \ .
  \end{align}
 Which in turn is equivalent to the existence of a line bundle $\mathcal{A}$ on $S$ such that
\begin{align}
{p_1}_*L\cong \underbrace{\bigoplus \mathcal{A}}_{N+1}\ .
 \end{align}

Below $\mbox{{Chow}}(\mathbf{X}\big / S)$ denotes the Chow form of the family $\mathbf{X}\big / S$,
 $\mu$ is the coefficient of $k^{n-1}$ in  $P(k)$, and
  $\mathcal{M}_{n}$ is the coefficient of $\binom{m}{n}$ in the CGKM expansion of
 $\det({p_1}_*L^{\otimes m})$ for $m>>0$ . A complete discussion of these notions is given in  ``\emph{CM Stability and the Generalized Futaki Invariant I}" . We refer the reader to that paper for the basic definitions and constructions that are used in the present article.

  We define an invertible sheaf on $S$ as follows.
  \begin{definition} (The Refined CM polarization \footnote{We use this terminology in order to distinguish this sheaf from one introduced by the second author in (\cite{psc}).}) \\
  \begin{align}\label{defn}
 {\mathbb{L}}_{1}(\mathbf{X}\big / S):= \{\mbox{\emph{Chow}}(\mathbf{X}\big / S)\otimes \mathcal{A}^{d(n+1)}\}^{n(n+1)+\mu}\otimes\mathcal{M}_{n}^{-2(n+1)}
 \end{align}
   \end{definition}
With the family $p_1: \mathbf{X}\rightarrow S$ fixed throughout, we will denote ${\mathbb{L}}_{1}(\mathbf{X}\big / S)$ by ${\mathbb{L}}_{1}$ in the remainder of the paper.

Our first result exhibits the Mabuchi energy as a \emph{singular} Hermitian metric on
${\mathbb{L}}_{1}$ .
\begin{theorem}
 Let $|| \ ||$ be any smooth Hermitian metric on $\mathbb{L}^{-1}_{1}$ .\footnote{$\mathbb{L}^{-1}_{1}$ denotes the dual of $\mathbb{L}_{1}$.}
 Then there is a continuous function $\Psi_{S}:S\setminus \Delta\rightarrow (-\infty, \ c)$
 such that for all $z\in S\big / \Delta$
 \begin{align}\label{singnrm}
&d(n+1)\nu_{\omega|_{\mathbf{X}_z}}(\vp_\sigma)= {\log}\left(e^{(n+1)\Psi_{S}(\sigma z)}\frac{||\ ||^{2}(\sigma z)}{||\ ||^{2}( z)}\right)\ .
\end{align}
Here $c$ denotes a constant which depends only on the choice of background K\"ahler metrics on $S$ and $\mathbf{X}$,  $\Delta$ denotes the discriminant locus of the map $p_1$, and $\omega |_{\mathbf{X}_z}$ denotes the restriction of the Fubini Study form of $\cpn$ to the fiber $\mathbf{X}_z$ .
 \end{theorem}
\begin{remark}
This should be compared with the main result in section 8 of \cite{psc}. The principal contribution of our present work is the observation that the whole theory in section 8 of \cite{psc} should be recast from the beginning with the sheaf  ${{\mathbb{L}}_{1}}$ .
\end{remark}
Let $X\hookrightarrow \cpn$ be an $n$ dimensional projective variety with Hilbert polynomial $P$. Let $Hilb_m(X)$ denote the $mth$ Hilbert point of $X$ (see \cite{sopv} for further information ). If $\lambda$ is a one parameter subgroup of $G$ then it is known (see \cite{sopv} ) that the weight, $w_{\lambda}(m)$,  of $Hilb_m(X)$ with respect to $\lambda$ is a \emph{polynomial} in $m$ of degree at most $n+1$. That is,
\begin{align*}
w_{\lambda}(m)=a_{n+1}(\lambda)m^{n+1}+a_n(\lambda)m^n+\dots \ .
\end{align*}
 Then the ratio may be expanded as follows.
 \begin{align*}
&\frac{w_{\lambda}(m)}{mP(m)} =F_{0}(\lambda)+F_{1}(\lambda)\frac{1}{m}+\dots +F_{l}(\lambda)\frac{1}{m^{l}}+\dots\\
 \end{align*}
\begin{definition}(Donaldson (\cite{skdtoric}))\newline
$F_{1}(\lambda)$ is the generalized Futaki invariant of $X$ with respect to $\lambda$ .
\end{definition}
 In our previous paper we have shown the following.\\
\ \\
 \textbf{Theorem} (\emph{The weight of the Refined CM polarization})\\
\ \\
i) \emph{There is a natural $G$ linearization on the line bundle ${{\mathbb{L}}_{1}}$} . \\
\ \\
ii) \emph{Let $\lambda$ be a one parameter subgroup of $G$. Let $z\in  \mathfrak{Hilb}_{\mathbb{P}^{N}}^{P}(\mathbb{C})$. Let $w_{\lambda}(z)$ denote the weight  of the restricted $\mathbb{C}^{*}$ action (whose existence is asserted in i)) on
${{\mathbb{L}}^{-1}_{1}}|_{z_{0}}$  where  $z_{0}=\lambda(0)z$. Then}
\begin{align}\label{weight}
{w_{\lambda}(z)=F_{1}(\lambda)}\ .
 \end{align}

 The main result of the paper is  the following corollary of (\ref{singnrm}) and (\ref{weight}) .
 \begin{corollary}(Algebraic asymptotics of the Mabuchi energy)\\
 Let $\vplt$ be the Bergman potential associated to an algebraic 1psg $\lambda$ of $G$, and let $z\in S\setminus\Delta$. Then there is an asymptotic expansion
\begin{align}\label{asymp}
d(n+1){ \nu_{\omega|_{\mathbf{X}_z}}(\varphi_{\lambda(t)})-\Psi_{S}({\lambda(t)})=
 F_{1}(\lambda)\log(|t|^2)+O(1)} \ \mbox{as}\ |t| \rightarrow 0. \quad
\end{align}
Moreover  $\Psi_{S}(\lambda(t))= \psi(\lambda)\log(|t|^{2})+O(1)$  where  $\psi({\lambda}) \in \mathbb{Q}_{\geq 0} $. Moreover,   $\psi({\lambda}) \in \mathbb{Q}_{+} $ if and only if ${\lambda(0)}\mathbf{X}_z=\mathbf{X}_{\lambda(0)z}$ (the limit cycle  \footnote {See \cite{sopv} pg. 61.} of $\mathbf{X}_z$ under $\lambda$ ) has a component of multiplicity greater than one. Here $O(1)$ denotes any quantity which is bounded as $|t|\rightarrow 0$.
\end{corollary}
Moser iteration and a refined Sobolev inequality (see \cite{simon} and \cite{sobolev}) yield the following.
\begin{corollary}
If $\nu_{\omega|_{\mathbf{X}_z}}  $ is proper (bounded from below) then the generalized Futaki invariant of $\mathbf{X}_z$ is {strictly negative} (nonnegative) for all $\lambda \in G$.   \end{corollary}

 \begin{remark}
We call the left hand side of (\ref{asymp}) the \emph{reduced} K-Energy along $\lambda$.
We also point out that while it is certainly the case that  $F_{1}(\lambda)$ may be defined for any subscheme of $\cpn$ it evidently only controls the behavior of the K-Energy when ${\lambda(0)\mathbf{X}_z} $ is  reduced.
\end{remark}

\begin{remark}
 {The precise constant $d(n+1)$ in front of $\nu_{\omega}$ is not really crucial, since what really matters  is the} sign \emph{of $F_{1}(\lambda)+ \psi({\lambda})$.}
 {That $\Psi_{S}(\lambda(t))$ has logarithmic singularities can be deduced from \cite{ags}.}
\end{remark}

\begin{remark}
 {We emphasize that we {do not} assume the limit cycle is smooth.}
\end{remark}
\begin{center}{\S 2  Background and Motivation}\end{center}
Let  $(X,\omega)$ be a compact K\"ahler manifold ($\omega$ not necessarily a Hodge class) and $P(X,\omega):= \{\varphi \in  C^{\infty}(X):\omega_{\varphi}:= \omega + \frac{\sqrt{-1}}{2\pi}\dl\dlb\varphi >0\}$  the space of K\"ahler potentials. This is the usual description of all K\"ahler metrics in the same class as $\omega$ (up to translations by constants).  It is not an overstatement to say that the most basic problem in K\"ahler geometry is the following
\begin{center}
\emph{Does there exist  $\varphi \in P(X,\omega)$ such that } $\mbox{Scal}(\omega_{\varphi})\equiv \mu$ . $(*)$\end{center}
This is a fully nonlinear \emph{fourth order} elliptic partial differential equation for $\varphi$.
  $\mu$ is a constant, the average of the scalar curvature, it depends only on $c_{1}(X)$ and $[\omega]$. When $c_{1}(X)>0$ and $\omega$ represents the \emph{anticanonical} class
 a simple application of the Hodge Theory shows that $(*)$ is equivalent to the \emph{Monge-Ampere equation}.
  \begin{align*}
  \frac{\mbox{det}(g_{i\overline{j}}+\varphi_{i\overline{j}})}{\mbox{det}(g_{i\overline{j}})}=e^{F-\kappa\varphi} \quad ( \kappa=1) \qquad (**)
  \end{align*}
 $F$ denotes the Ricci potential. When $\kappa=0$ this is the celebrated Calabi problem solved by S.T.Yau
 and when $\kappa < 0$ this was solved by Aubin and Yau independently in the 70's.
 It is well known that $(*)$ is actually a \emph{variational} problem. There is a natural energy  on the space $P(X,\omega)$ whose critical points are those $\varphi$ such that $\omega   _{\varphi}$ has constant scalar curvature (csc).
This energy was introduced by T. Mabuchi (\cite{mabuchi}) in the 1980's. It is called the \emph{\textbf{K-Energy map}} (denoted by $\nu_{\omega}$) and is given by the following formula
\begin{align*}
 \qquad \nu_{\omega}(\varphi):= -\frac{1}{V}\int_{0}^{1}\int_{X}\dot{\varphi_{t}}(\mbox{Scal}(\varphi_{t})-\mu)\omega_{t}^{n}dt.
\end{align*}
Above, $\varphi_{t}$ is a smooth path in $ P(X,\omega)$ joining $0$ with $\varphi$. The K-Energy does not depend on the path chosen.
In fact there is the following well known formula for $\nu_{\omega}$ where $O(1)$ denotes a quantity which is bounded on $P(X,\omega)$.
\begin{align*}
&\nu_{\omega}(\varphi)=\int_{X}\mbox{log}\left(\frac{{\omega_{\varphi}}^{n}}{\omega^{n}}\right)\frac{{\omega_{\varphi}}^{n}}{V} - \mu(I_{\omega}(\varphi)-J_{\omega}(\varphi))+O(1) \\
& J_{\omega}(\varphi):= \frac{1}{V}\int_{X}\sum_{i=0}^{n-1}\frac{\sqrt{-1}}{2\pi}\frac{i+1}{n+1}\dl\varphi \wedge \dlb
\varphi\wedge \omega^{i}\wedge {\omega_{\varphi} }^{n-i-1}\\
&I_{\omega}(\varphi):= \frac{1}{V}\int_{X}\varphi(\omega^{n}-{\omega_{\varphi}}^{n})
\end{align*}

 We have written down the K-energy in the case when $\omega = c_{1}(X)$ . Observe that $\nu_{\omega}$ is essentially the \emph{difference} of two positive terms. What is of interest for us is that  the problem $(*)$ is not only a {variational} problem but a \emph{minimization} problem. With this said we have the following fundamental result. \\
\ \\
 \noindent \textbf{Theorem} (S. Bando and T. Mabuchi \cite{bandmab}) \\
 \emph{If $\omega=c_{1}(X)$ admits a K\"ahler Einstein metric then $\nu_{\omega}\geq 0$. The absolute minimum is taken on the solution to $(**)$ (which is unique up to automorphisms of $X$).}\\
\ \\
Therefore a \emph{necessary} condition for the existence of a K\"ahler Einstien metric is a bound from below on $\nu_{\omega}$. In order to get a \emph{sufficient} condition
one requires that the K-energy \emph{grow} at a certain rate. Precisely, it is required that the K-Energy be \emph{proper}.
This concept was  introduced by the second author  in \cite{psc}.
\begin{definition}
{$\nu_{\omega}$ is \textbf{proper} if there exists a strictly increasing function $f:\mathbb{R}_{+}\ra\mathbb{R}_{+}$ (where $\lim_{T\ra \infty}f(T)=\infty$) such that $\nu_{\omega}(\varphi)\geq f(J_{\omega}(\varphi))$ for all $\varphi\in P(M,\omega)$}.
 \end{definition}

\noindent \textbf{Theorem} (\cite{psc})\\
  \emph{ Assume that $Aut(X)$ is discrete. Then $\omega=c_{1}(X)$ admits a K\"ahler Einstein metric if and only if $\nu_{\omega}$ is proper}.\\
\ \\
 The next result was established by the second author and Xiuxiong Chen. It holds in an \emph{arbitrary} K\"ahler class $\omega$. An alternative proof of this was given by Donaldson for polarized projective
 manifolds.\\
\ \\
 \textbf{Theorem} (\cite{chentian})\\
  \emph{ If $\omega$ admits a metric of csc then} $\nu_{\omega}\geq 0$.\\

 In this paper our interest is to test for a lower bound of $\nu_{\omega}$ along the large but finite dimensional group $G$  of \emph{matrices} in the polarized case.  When we restrict our attention to $G$ we make the connection with Mumfords' Geometric Invariant Theory. The past couple of years have witnessed quite a bit of activity on this problem due to this connection.

To put things in historical perspective consider the various formulations of the Futaki invariant.\\
\ \\
i) 1983 Futaki (\cite{futaki}) introduces his invariant as a lie algebra character on a Fano manifold $X$
\begin{align*}
F_{\omega}:\eta(X)\ra \mathbb{C}.
\end{align*}

\noindent ii) 1986 Mabuchi (see \cite{mabuchi} ) integrates the Futaki invariant with the introduction of the K-energy map. The linearization of the K-energy along orbits of holomorphic vector fields is the real part of the Futaki invariant.\\
\ \\
iii) 1992  Ding and Tian (\cite{dingtian}) introduced the \emph{generalized} Futaki invariant. Here the \emph{jumping of complex structures} is introduced. The limit of the derivative of the K-Energy map is identified with the generalized Futaki invariant of  $X^{\lambda(0)}$ provided this limit has at most \emph{normal} singularities.\\
\ \\
iv) 1997  The CM polarization is defined (see \cite{psc}) for \emph{smooth} families, as the relative canonical bundle is explicitly involved in the definition. K-Stability is defined in terms of special degenerations and the generalized Futaki invariant. \\
\ \\
v) 1999 Yotov formulated the generalized Futaki Invariant in terms of equivariant Chow groups of a \emph{normal} variety.\\
\ \\
vi) 2002 For an \emph{arbitrary scheme} Donaldson (\cite{skdtoric}) defined the weight $F_{1}(\lambda)$. This is identified with the limit of the derivative of K-energy (by \cite{dingtian}) when the limit cycle is  a \emph{smooth} (or normal) scheme.\\

\begin{remark}
 We hope that we have clarified the role of the CM polarization. The main point is that once the CM polarization is extended to the Hilbert scheme (\cite{cms1}) the polarization computes the precise asymptotics of the K-energy of any generic fiber of the map $\mathbf{X}\rightarrow S$. This extension was made possible by an application of the Knudsen Mumford expansion of the determinant of direct images of perfect complexes of sheaves (see \cite{detdiv}). In fact, $\psi(\lambda)$ already appeared in work of the second author (see \cite{psc}). Despite this, the role of $\psi(\lambda)$ becomes more precise 
 in the present work. 
\end{remark}
\section{Algebraic potentials}
In order to connect these notions to the K-Energy map we now  give an account of how to associate an admissible potential $\varphi_{\lambda(t)}$ to a one parameter subgroup of $G$.
In order to detect properness (conjecturally) one restricts attention to the subspace of \emph{Bergman metrics} inside $P(M,\omega)$ since these metrics are \emph{dense} in $P(M,\omega)$ (see \cite{tianberg}, \cite{ruan}, \cite{zel}, \cite{cat}).  By definition these metrics are induced by the Kodaira embeddings furnished by the polarization $L$. The construction is as follows.
 We have an embedding
\[
\begin{CD}
X@>>{ L^{}}>\mathbb{P}(H^{0}(X, L^{})^{*})= \mathbb{P}^{N}
\end{CD}
\]
furnished by some basis $\{S_{0},\dots,S_{N}\}$ of $H^{0}(X, L)$.
Observe that with the natural Hermitian metric on $H^{0}(X, L)$, the induced Fubini-Study metric on $ \mathbb{P}^{N}$ is related to the curvature of the initial metric on $ L$ by the formula
\begin{align*}
\omega_{FS}|_{X}=\omega +\frac{\sqrt{-1}}{2\pi}\dl\dlb \mbox{log}\left(\sum_{i=0}^{N}||S_{i}||^{2}\right).
\end{align*}
We conclude that
\begin{align*}
\mbox{log}\left(\sum_{i=0}^{N}||S_{i}||^{2}\right)\in P(X,\omega).
\end{align*}
Let $\sigma \in SL(N+1,\mathbb{C})$, then
\begin{align*}
\sigma^{*}(\omega_{FS})=\omega_{FS} +\frac{\sqrt{-1}}{2\pi}\dl\dlb\varphi_{\sigma}  .
\end{align*}
Where $\varphi_{\sigma}$ is given by the formula
\begin{align*}
\varphi_{\sigma}=\mbox{log}\left(\frac{||\sigma z||^{2}}{|| z||^{2}}\right).
\end{align*}
We let $\{T_{0},\dots T_{N}\}$ denote the corresponding change of basis
\begin{align*}
\begin{pmatrix}\sigma_{00} & \dots & \sigma_{0N}\\
        \sigma_{10} & \dots & \sigma_{1N}  \\
         \dots& \dots & \dots \\
         \sigma_{N0} & \dots & \sigma_{NN}   \\ \end{pmatrix}
\begin{pmatrix}S_{0}\\
               \dots \\
               \dots \\
S_{N}\\ \end{pmatrix}
=
\begin{pmatrix}T_{0}\\
               \dots \\
               \dots \\
T_{N}\\ \end{pmatrix} \ .
\end{align*}
Then we have
\begin{align*}
{\varphi_{\sigma}}|_{X}=
\mbox{log}\left(\frac{\sum_{i=0}^{N}||T_{i}||^{2}}{\sum_{i=0}^{N}||S_{i}||^2}\right)\ .
\end{align*}
Putting these ingredients together gives
\begin{align}\label{bergman}
\sigma^{*}\omega_{FS}|_{X}=\omega +\frac{\sqrt{-1}}{2\pi}\dl\dlb\mbox{log}\left({\sum_{i=0}^{N}||T_{i}||^{2}}\right)\ .
\end{align}
Therefore, if we \emph{fix} a basis of $H^{0}(X, L)$ we get a natural map
\begin{align*}
SL(N+1,\mathbb{C})\rightarrow P(X,\omega)\ .
\end{align*}
A \emph{one parameter subgroup} of $SL(N+1,\mathbb{C})$ is an algebraic\footnote{``algebraic'' means that the matrix coefficients $\lambda(t)_{i,j}\in \mathbb{C}[t,t^{-1}]$.} homomorphism
\begin{align*}
\lambda:\mathbb{C}^{*}\rightarrow SL(N+1,\mathbb{C})\ .
\end{align*}
 Any such $\lambda(t)$ can be diagonalised. That is, we may assume that $\lambda(t)$ takes values in the standard maximal torus $H\cong (\mathbb{C}^*)^N$ of $\slnc$.
\begin{align*}
\lambda(t)=\begin{pmatrix}t^{m_{0}}&\dots&\dots& 0\\
                             0&t^{m_{1}}&\dots& 0\\
                             0&\dots&\dots& t^{m_{N}}
                             \end{pmatrix}\ .
\end{align*}
The exponents $m_{i}$ satisfy
\begin{align*}
 \quad \sum_{0\leq i \leq N}m_{i}=0.
\end{align*}
We arrive at the following formula.
 \begin{align*}
\varphi_{\lambda(t)}(z):= \frac{1}{k}\mbox{log}\left(\sum_{0\leq j\leq N}|t|^{2m_{j}}||S_{j}||^{2}(z)\right).
\end{align*}
 Now we may consider the K-energy map as a function on $\slnc$ .
\section{Singular Hermitian metrics}
 \subsection{Proof of Theorem 1} In part I of this work the authors provided the following  formula for the first Chern class of ${\mathbb{L}}_{1}$ .
 \smallskip
\begin{align}\label{c1}
c_{1}( {\mathbb{L}}_{1})={p_1}_*\left((n+1)c_1(K_{\mathbf{X}/ S})c_1({L})^n+\mu\ c_1({L})^{n+1}\right) \quad K_{\mathbf{X} / S}:= K_{\mathbf{X}}\otimes {p^*_1}(K^{\vee}_S)\ .
\end{align}
( \ref{c1}) allows us to exhibit the K-energy map as a \emph{singular} metric on the CM polarization (see \cite{psc}).
 Recall that $p^{-1}(z)=\mathbf{X}_{z}\subset \cpn$, where $z\in {S}_{\infty}:=S\setminus \Delta$.
We define
\begin{align*}
G\mathbf{X}_z:= \{(\sigma,y)\in G \times \cpn:y\in \sigma \mathbf{X}_z\} \ .
\end{align*}
Observe that $G\mathbf{X}_{z}$ is biholomorphic to $G\times \mathbf{X}_{z}$.
Then we have the following diagram, where $p_{z}$ denotes the evaluation map, i.e.
$p_{z}(\sigma):= \sigma z$.

 \begin{diagram}
p_{z}^{*}(\mathbf{X})\cong G\mathbf{X}_z&\rTo^{p_{z,2}}&\mathbf{X}&\rInto^{\iota}& S\times \cpn&\rTo^{p_{2}}&\cpn \\
\dTo ^{p_{z,1}}&&\dTo^{p}&\ldTo_{\pi}&\\
  G&\rTo ^{p_{z}}&S&
\end{diagram}
Given $z\in B\setminus \Delta$ we can consider $K_{\mathbf{X}_z}$, the canonical bundle of the fiber $\mathbf{X}_z$. These fit together holomorphically into a line bundle $K_{\infty}$ on $\mathbf{X}\setminus p^{-1}(\Delta)$.
On the other hand, the relative canonical bundle  $K_{p}$ of the map $p$  exists and lives on \emph{all} of $\mathbf{X}$.
\begin{align*}
{K}_{p}:= K_{\mathbf{X}}\otimes p^{*} K_{S}^{-1}
\end{align*}
When we restrict this sheaf to $\mathbf{X}\setminus p^{-1}(\Delta)$ we have an isomorphism
\begin{align*}
{K}_{p}  \cong K_{\infty}\ .
\end{align*}

  $\iota^{*}p_{2}^{*}\omega_{FS}$ restricts to a K\"ahler  metric on $p^{-1}(z)$ $(z \in {S}_{\infty})$ and hence induces a Hermitian metric on the bundle $K_{\infty}$.  We denote
 its curvature
 by $R( \iota^{*}p_{2}^{*} (\omega_{FS}))$.
Let $g_{\mathbf{X}}$ and $g_S$ denote two K\"ahler metrics on
$\mathbf{X}$ and $S$ respectively.
In this way we obtain a metric on the relative canonical bundle $K_{p}$.
We let $R_{f}$ denote its curvature
\begin{align*}
R_{p}:= R(g_{\mathbf{X}})- p^{*}R(g_S).
\end{align*}
 In this way we obtain \emph{two} metrics on the relative canonical bundle over the smooth locus.
 The crucial point is the following fact.\\
 \ \\
\begin{center} \emph{The curvatures of these metrics are {not} the same .}\end{center}
\ \\
The relation between them is given in the following proposition.
\begin{proposition} (`` {$\dl\dlb$ lemma along the fibers}'')\\
 {There is a smooth function $\Psi:\mathbf{X}\setminus p^{-1}(\Delta)\rightarrow \mathbb{R}$ such that}
\begin{align*}
\begin{split}
& 1) \ R(g_{\mathbf{X}})-p^{*}R(g_{ S}) + \frac{\sqrt{-1}}{2\pi}\dl\dlb \Psi = R( \iota^{*}p_{2}^{*} (\omega_{FS}))\ . \\
\ \\
& 2) \ \Psi \leq C,\ \mbox{for some constant $C$}\ .
\end{split}
\end{align*}
\end{proposition}
\begin{example} ({The universal family of hypersurfaces of degree $d$ in $\mathbb{C}P^{n+1}$}) \\
 \begin{align*}
 &S:=\mathbb{P}(H^0(\mathbb{C}P^{n+1},\mathcal{O}(d)))\\
 \ \\
 &\mathbf{X}:=\{([f],[z])\in S\times \mathbb{C}P^{n+1}\ | f(z)=0\}\\
 \ \\
 &p:= p_1 \qquad (\mbox{projection onto the first factor}) \ .
 \end{align*}
 Let $ |||\ . |||$ denote any norm on $H^0(\mathbb{C}P^{n+1},\mathcal{O}(d))$, with associated Fubini-Study metric $\omega_{S}$. Then a computation shows that
 \begin{align*}
\Psi(([f],[z]))= \log\left(\frac{\sum_{i=0}^{n+1}|\frac{\dl f}{\dl z_i}(z)|^2}{|||f|||^2||z||^{2(d-1)}}\right)\ .
\end{align*}
\end{example}
The next result is a \emph{pointwise} version of  (\ref{c1}) .
 \begin{proposition}
 {There is a continuous Hermitian metric $||\ ||$ on $\mathbb{L}^{-1}_1$ such that, in the sense of currents we have}
\begin{align*}
\frac{\sqrt{-1}}{2\pi}\dl\dlb\emph{log}({||\ ||^{2}})= (n+1)p_{*}(R(g_{\mathbf{X}})-p^{*}R(g_{ S}))p_{2}^{*}(\omega_{FS})^{n} + \mu p_{*}p_{2}^{*}(\omega_{FS})^{n+1}\ .
\end{align*}
\end{proposition}
\begin{proof}
See Proposition 4.3 pg. 2576 of \cite{ags}  .
\end{proof}
Now we pull back the curvature form of $K_{\infty}$ to $G\mathbf{X}_z$
 \begin{align*}
R_{G|\mathbf{X}_z}:= p_{z,2}^{*}(R(\pi_{2}^{*}(\omega_{FS}))).
\end{align*}
Recall that for $\sigma \in G$ we define $\varphi_{\sigma}$ by the relation
\begin{align*}
\sigma^{*}\omega_{FS}= \omega_{FS}+\frac{\sqrt{-1}}{2\pi}\dl\dlb \varphi_{\sigma}.
\end{align*}
 Let  $\nu_{\omega,z}(\sigma)$ denote the K energy of $(\mathbf{X}_z, \omega_{FS})$ applied to the potential $\varphi_{\sigma}$.
With these notations in place we have the following result.\\
 \begin{proposition}\label{ddb} (The complex Hessian of the K-Energy map on $G$)\\
 {For every  smooth compactly supported  $(N^{2}+2N-1,N^{2}+2N-1)$ form $\eta$ on $G$ we have}
\begin{align*}
d(n+1)\int_{G}\nu_{\omega,z}(\varphi_{\sigma} )\dl\dlb \eta = \int_{G\mathbf{X}_z}((n+1)R_{G|\mathbf{X}_z}+{\mu}p_{2}^{*}(\omega_{FS}))\wedge p_{2}^{*}(\omega_{FS})^{n}\wedge p_{z,1}^{*}\eta.
\end{align*}
 \end{proposition}
  The proof of proposition 3 appears in the next section after some standard preliminaries on Bott Chern classes.
\subsection{Bott Chern secondary classes}
Let $\phi$ be a $GL_N(\mathbb{C})$ invariant polynomial on $M_{N\times N}(\mathbb{C})$ homogeneous of degree $d$. $\phi_1$ denotes the \emph{complete polarization} of $\phi$.
Let $E$ be a holomorphic vector bundle
of rank $N$ over a base $X$. Let $h_1$ and $h_0$ be two Hermitian
metrics on $E$ and $\frac{\sqrt{-1}}{2\pi}R(h_i)$ the curvatures.
Then we define the \emph{Bott-Chern class} $BC(\phi,E;h_0,h_1)$ by the expression
 \begin{align}
\label{goog}
BC(\phi,E;h_0,h_1):= \int_0^1\phi_1(h_{t}^{-1}\dot h_t,\overbrace{\frac{\sqrt{-1}}{2\pi}R_t,\dots,\frac{\sqrt{-1}}{2\pi}R_t}^{d-1})dt \in \Omega_X^{(d-1,d-1)}
\end{align}
    $h_t$ is any piecewise
$C^{1}$ path of Hermitian metrics joining $h_0$ and $h_1$.
The point of the construction is the following identity.
\begin{align*}
\frac{\sqrt{-1}}{2\pi}\dl\dlb BC(\phi,E;h_0,h_1)= \left(\frac{\sqrt{-1}}{2\pi}\right)^d(\phi(R_{h_0})-\phi(R_{h_1}))
\end{align*}
Let $d=n+1$ where, $n=dim(X)$ in this case $BC(\phi,E;h_0,h_1)$ has top dimension and we may introduce the \emph{Donaldson Functional associated to $\phi$} .
\begin{align}
D_E(h_0,h_1):=\int_{X}BC(\phi,E;h_0,h_1)
\end{align}
When $h_0$ is fixed, we consider it to be a functional on
$\mathcal{M}_{E}$ (the space of hermitian metrics on $E$). In what
follows we take $\phi = Ch_{n+1}$, the $n+1^{st}$ component of the
chern character. We can extend the Donaldson functional to
``virtual bundles'' $\mathcal{E}=E-F$ by observing that
 a Hermitian metric $h$ on $\mathcal{E}$ is just a pair of metrics, one on $E$ and one on $F$.
\begin{align*}
h=(h^{E},h^{F}) \ .
 \end{align*}
 We set
 \begin{align}
 BC(\phi,\mathcal{E};h_0,h_1):= BC(\phi,E;h^{E}_0,h^{E}_1)- BC(\phi,F;h^{F}_0,h^{F}_1)
\end{align}
Let $h:Y\rightarrow \mathcal{M}_{\mathcal{E}}$ be a smooth map, where $Y$ is a complex manifold of dimension $m$.
\begin{lemma} Let $\phi$ be homogeneous of degree $n+1$ and $h_0$ a fixed metric on $\mathcal{E}$. Then for all smooth compactly supported forms $\psi$ of type $(m-1,m-1)$ we have the identity
\begin{align}\label{workhorse}
\frac{\sqrt{-1}}{2\pi}\int_YD_\mathcal{E}(\phi; h_0, h(y))\dl_Y\dlb_Y\psi=\int_{Y\times X}\phi(R(\frac{\sqrt{-1}}{2\pi}h(y)))\wedge \pi_1^*(\psi)\ .
\end{align}
\end{lemma}
Next we want to realize the Mabuchi K-energy as the
Donaldson functional, with respect to the polynomial $\phi = Ch_{n+1}$, of a certain virtual bundle to be defined below. Then proposition (\ref{ddb}) follows at once from the preceding lemma.

Let $X$ be a complex projective manifold (in our present application $X$ is a smooth fiber of $\mathbf{X}\overset{p}{\rightarrow} S$), and let $L$ be
the restriction of $\underline{O}(1)$ to $X$. Let $\varphi$ be a
kahler potential. The two metrics $h_{FS}$ and
$e^{-\varphi}h_{FS}$ induce metrics on the
canonical bundle $\mathcal{K}$. We consider the virtual
bundle
\begin{align}
2^{n+1}\mathcal{E}:=
(n+1)(\mathcal{K}^{-1}-\mathcal{K})(L-L^{-1})^{n}-\mu(L-L^{-1})^{n+1} \ .
\end{align}
Here $\mu$ is the average of the scalar curvature. We need to
calculate the following terms.
 \begin{align}
 \begin{split}
&BC(\phi; \mathcal{K}^{-1}\otimes L^{n-2j},h_0,h_1) \\
\ \\
&BC(\phi; \mathcal{K}\otimes L^{n-2j},h_0,h_1) \\
 \ \\
 &BC(\phi;L^{n+1-2j},h_{0},h_{1}) \ .
 \end{split}
 \end{align}
The path of metrics for the first two expression are given as follows.
\begin{align}
\begin{split}
& h_{\mathcal{K}^{-1}\otimes L^{n-2j},t}:= \det(g_{\alpha
\overline{\beta}}+t\frac{\dl^{2}}{\dl z_{\alpha}\dl
\overline{z_{\beta}}}\varphi)e^{-t(n-2j)\varphi}h_{FS}^{n-2j} \\
\ \\
& h_{\mathcal{K}\otimes L^{n-2j},t}:= \det(g_{\alpha
\overline{\beta}}+t\frac{\dl^{2}}{\dl z_{\alpha}\dl
\overline{z_{\beta}}}\varphi)^{-1}e^{-t(n-2j)\varphi}h_{FS}^{n-2j} \ .
\end{split}
\end{align}
The complete polarization of $\phi$ is given by
\begin{align}
\phi_{1}(B,A \dots A)= tr(BA^{n}) \quad A,B \in M_{k}(\mathbb{C}) \ .
\end{align}
Therefore ,
\begin{align}
\begin{split}
BC(\mathcal{K}^{-1}\otimes
L^{n-2j},h_0,h_1)=\int_{0}^{1}(\Delta_{t\varphi}\varphi -
(n-2j)\varphi)((n-2j)\omega_{t\varphi} + Ric_{\omega_{t}})^ndt \\
\ \\
 BC(\mathcal{K}\otimes L^{n-2j},h_0,h_1)=
-\int_{0}^{1}( \Delta_{t\varphi}\varphi +
(n-2j)\varphi)((n-2j)\omega_{t\varphi} - Ric_{\omega_{t}})^ndt \ .
\end{split}
\end{align}
 Similarly we have
 \begin{align}
BC(L^{n+1-2j},h_{0},h_{1})= -(n+1-2j)^{n+1}\int_{0}^{1}\varphi\omega_{t}^{n}dt \quad \omega_{t}:= \omega+t\dl \dlb \varphi \ .
\end{align}

We see that
\begin{align}
BC((L-L^{-1})^{n+1},h_{FS},e^{-\varphi}h_{FS})=
-\sum_{j=0}^{n+1}(-1)^{j}\binom{n+1}{j}(n+1-2j)^{n+1}\int_{0}^{1}\varphi\omega_{t}^{n}dt  \ .
\end{align}
Now we need the following numerical identity.
\begin{align}\label{partout}
\sum_{j=0}^{n+1}(-1)^{j}\binom{n+1}{j}(n+1-2j)^{i}=\begin{cases}
0 & i<n+1\  \mbox{or} \ i=n+2\\
(n+1)!2^{n+1} & i= n+1\ . \\
\end{cases}
\end{align}
It follows at once that
\begin{align}\label{aub}
\int_{X}BC((L-L^{-1})^{n+1},h_{FS},e^{-\varphi}h_{FS})=-(n+1)!2^{n+1}\int_{0}^{1}\int_X\varphi\omega_{t}^{n}dt \ .
\end{align}
It follows from (3.9) that
\begin{align*}
BC(\mathcal{K}^{-1}\otimes L^{n-2j})= \int_{0}^{1}\Delta_{t\varphi}\varphi\sum_{i=0}^{n}\binom{n}{i}(n-2j)^{i}Ric_{t}^{n-i}\omega_{t}^{i}-\int_{0}^{1}\sum_{i=0}^{n}\binom{n}{i}(n-2j)^{i+1}\varphi Ric_{t}^{n-i}\omega_{t}^{i}\ .
\end{align*}
We use the identity (\ref{partout}) to see that
\begin{align}\label{antican}
\sum_{j=0}^{n}(-1)^{j}\binom{n}{j}BC(\mathcal{K}^{-1}\otimes L^{n-2j})=n!2^{n}\int_{0}^{1}\left(\Delta_{t}\varphi \omega_{t}^{n}-\varphi nRic_{t}\omega_{t}^{n-1}\right)dt \ .
\end{align}
Similarly we have the second term
\begin{align}\label{can}
\sum_{j=0}^{n}(-1)^{j+1}\binom{n}{j}BC(\mathcal{K}\otimes L^{n-2j})= n!2^{n} \int_{0}^{1}\left(\Delta_{t}\varphi \omega_{t}^{n}-\varphi nRic_{t}\omega_{t}^{n-1}\right)dt \ .
\end{align}
The next lemma follows at once from summing up (\ref{antican}), (\ref{can}), and (\ref{aub}) .
\begin{lemma}\label{dnfnctl}
Let $D(\mathcal{E},h_{FS},e^{-\varphi}h_{FS})$ denote the Donaldson functional of $Ch_{n+1}$ with respect to $\mathcal{E}$. Then the following identity holds.
\begin{align}
D(\mathcal{E},h_{FS},e^{-\varphi}h_{FS}) = \nu_{\omega}(\varphi)
\end{align}
\end{lemma}
Let $\varphi=\varphi_{\sigma}$ and apply \ref{workhorse} to lemma \ref{dnfnctl} to conclude the proof of proposition \ref{ddb} . $\Box$

Next we observe that the identity
 \begin{align}
R_{G|\mathbf{X}_z}=p_{2,z}^{*}\left(R(g_{\mathbf{X}}) -p^{*}R(g_{ S}) +
\frac{\sqrt{-1}}{2\pi}\dl\dlb \Psi \right)
\end{align}
 together with the previous lemmas yields the following corollary.
  \begin{corollary}
 {The function }
 \begin{align*}
\sigma \in G \rightarrow D(\sigma):= d(n+1)\nu_{\omega,z}(\sigma)  -{\log}\left(e^{(n+1)\Psi_{S}(\sigma z)}\frac{||\ ||^{2}(\sigma z)}{||\ ||^{2}( z)}\right)
\end{align*}
 {is pluriharmonic.  Where we have defined}  $\Psi_{S}(z):= \int_{\{y\in f^{-1}(z)\}}\Psi(y) p_{2}^{*}(\omega_{FS})^{n}$. \\
Moreover $\Psi_{S}(z)\leq C$ on $S\setminus \Delta$ , extends continuously to the locus of reduced fibers, and $\lim_{z\rightarrow z_{\infty}}\Psi_{S}(z)= -\infty$ whenever $\mathbf{X}_{z_{\infty}}$ is non-reduced.
\end{corollary}
\begin{remark}
The construction of $\Psi$ and $\Psi_S$ as well as their behavior on the locus of singular fibers can be seen directly in example 1. The general case is treated in lemma 8.5 pg. 31 in \cite{psc} .
\end{remark}
Since $\pi_{1}(G)=1$ there is a (nonvanishing) entire function $\xi$ on $G$ such that
\begin{align*}
D(\sigma)= \log(|\xi(\sigma)|^{2}).
\end{align*}
An analysis of the growth of  this function on the standard compactification $\overline{G}$
\begin{align*}
\overline{G}:= \{[(w_{ij},\ z)] \in \mathbb{P}^{(N+1)^{2}}: \mbox{det}(w_{ij})=z^{N+1} \}
\end{align*}
reveals that it must reduce to a \emph{constant}.

Tying everything together establishes our main result.

\textbf{Theorem 1}
(\emph{The K-Energy as a singular metric on $\mathbb{L}^{-1}_1$})
 \begin{align}\label{mainthm}
d(n+1)\nu_{\omega,z}(\sigma)= {\log}\left(e^{(n+1)\Psi_{S}(\sigma z)}\frac{||\ ||^{2}(\sigma z)}{||\ ||^{2}( z)}\right)\ .
\end{align}

 We proceed to the proof of corollary 1 . First substitute $\sigma=\lambda(t)$ in (\ref{mainthm}). Then we have the string of identities.
 \begin{align*}
d(n+1)\nu_{\omega,z}(\lambda(t))&= {\log}\left(e^{(n+1)\Psi_{S}(\lambda(t) z)}\frac{||\ ||^{2}(\lambda(t)z)}{||\ ||^{2}( z)}\right)\\
\ \\
&= (n+1)\Psi_{S}(\lambda(t) z) + {\log}\left(\frac{||\ ||^{2}(\lambda(t)z)}{||\ ||^{2}( z)}\right)\\
\ \\
&= (n+1)\Psi_{S}(\lambda(t) z) + {\log}\left(\frac{||\ ||^{2}(t^{w_{\lambda}(z)-w_{\lambda}(z)}\lambda(t)z)}{||\ ||^{2}( z)}\right)\\
\ \\
&= (n+1)\Psi_{S}(\lambda(t) z) +w_{\lambda}(z)\log(|t|^2)+O(1)\\
\ \\
&=F_{1}(\lambda)\log(|t|^2)+ (n+1)\Psi_{S}(\lambda(t) z) +O(1)\ .
\end{align*}

The passage from line 3 to 4 follows from the defining property of the weight (see the introduction to \cite{cms1}). The passage from line 4 to 5 is the statement of (\ref{weight}).
This completes the proof of corollary 1. $\Box$
   \subsection{Properness Implies that $F_1(\lambda)<0$ }
  Let $X:= \mathbf{X}_z$ a smooth fiber of $p$. Recall that the algebraic potential associated to a one parameter subgroup $\lambda$ is given by
 \begin{align*}
 \varphi_{t}:=\varphi_{\lambda(t)}=\mbox{log}(\sum_{i=0}^{N}t^{2q_{i}}||S_{i}||^{2}).
 \end{align*}

  Then, as we have seen, $\varphi _{t}\in P(X,\omega)$.
 Following Yau \cite{cma}, our plan is to use the standard Moser iteration to control $Osc(\varphi _{t})$ by $I_{\omega}(\varphi _{t})$.
Define
\begin{align*}
\varphi_{-}:= \mbox{Max}\{-\varphi_{t},1\}\geq 1.
\end{align*}
Let $p\in \mathbb{Z}_{+}$. Then we have the (obvious ) inequality
\begin{align*}
\varphi_{-}^{p}\frac{\sqrt{-1}}{2\pi}\dl\dlb\varphi \wedge\omega_{\varphi}^{n-1}\leq \varphi_{-}^{p}\omega_{\varphi}^{n}.
\end{align*}
Trivially this implies
\begin{align*}
\int_{X}\varphi_{-}^{p}\frac{\sqrt{-1}}{2\pi}\dl\dlb\varphi \wedge\omega_{\varphi}^{n-1}\leq \int_{X}\varphi_{-}^{p}\omega_{\varphi}^{n}
\leq \int_{X}\varphi_{-}^{p+1}\omega_{\varphi}^{n}.
\end{align*}
Next integrate by parts on the leftmost side of this inequality
\begin{align*}
\begin{split}
\int_{X}\varphi_{-}^{p}\frac{\sqrt{-1}}{2\pi}\dl\dlb\varphi \wedge\omega_{\varphi}^{n-1}&=-\int_{X}\frac{\sqrt{-1}}{2\pi}\dl\varphi_{-}^{p}\wedge\dlb\varphi\ \omega_{\varphi}^{n-1}\\
&=\int_{X}\frac{\sqrt{-1}}{2\pi}\dl\varphi_{-}^{p}\wedge\dlb\varphi_{-}\omega_{\varphi}^{n-1}\\
&=\frac{4p}{(p+1)^{2}}\frac{\sqrt{-1}}{2\pi}\int_{X}\dl\varphi_{-}^{\frac{p+1}{2}}\wedge\dlb\varphi_{-}^{\frac{p+1}{2}}\wedge\omega_{\varphi}^{n-1}\ .
\end{split}
\end{align*}
Since $\varphi_{-}\geq 1$ we deduce the gradient estimate
\begin{align*}
\frac{4p}{(p+1)^{2}}\frac{\sqrt{-1}}{2\pi}\int_{X}\dl\varphi_{-}^{\frac{p+1}{2}}\wedge\dlb\varphi_{-}^{\frac{p+1}{2}}\wedge\omega_{\varphi}^{n-1}
\leq \int_{X}\varphi_{-}^{p}\omega_{\varphi}^{n}\leq \int_{X}\varphi_{-}^{p+1}\omega_{\varphi}^{n}\ .
\end{align*}
We concentrate on the outermost inequality
\begin{align*}
\frac{4p}{n(p+1)^{2}}\int_{X}||\nabla_{\varphi_{t}}\varphi_{-}^{\frac{p+1}{2}}||^{2}_{\varphi_{t}}\omega_{\varphi_{t}}^{n}\leq \int_{X}\varphi_{-}^{p+1}\omega_{\varphi}^{n}\ .
\end{align*}
Now we invoke the Sobolev inequality
\begin{align*}
\left(\int_{X}\varphi_{-}^{\frac{(p+1)n}{n-1}}\frac{\omega_{\varphi}^{n}}{V}\right)^{\frac{n-1}{n}}\leq \mathcal{C}(\varphi_{t})\left(\int_{X}||\nabla_{\varphi_{t}}\varphi_{-}^{\frac{p+1}{2}}||^{2}_{\varphi_{t}}\frac{\omega_{\varphi_{t}}^{n}}{V}+\int_{X}\varphi_{-}^{p+1}\frac{\omega_{\varphi}^{n}}{V}\right) \ .
\end{align*}
$\mathcal{C}(\varphi_{t})$ is the Sobolev constant of the metric $\omega +\frac{\sqrt{-1}}{2\pi}\dl\dlb\varphi_{t}$.
Concerning this constant we have the crucial
\begin{proposition}(\cite{simon}, \cite{sobolev})
There is a positive constant $\delta=\delta(n)$ such that for all $\sigma\in \slnc$ we have
\begin{align*}
\mathcal{C}(\varphi_{\sigma}) <\delta.
\end{align*}
\end{proposition}
This follows from the fact the complex projective subvarieties are \emph{minimal} as Riemannian submanifolds of $\cpn$ and hence have vanishing mean curvature.

Therefore,
\begin{align*}
\left(\int_{X}\varphi_{-}^{\frac{(p+1)n}{n-1}}\frac{\omega_{\varphi}^{n}}{V}\right)^{\frac{n-1}{n}}\leq n(p+1)\delta\int_{X}\varphi_{-}^{p+1}\frac{\omega_{\varphi}^{n}}{V}.
\end{align*}
Now extract the $p+1$st root of both sides to get
\begin{align*}
\left(\int_{X}\varphi_{-}^{\frac{(p+1)n}{n-1}}\frac{\omega_{\varphi}^{n}}{V}\right)^{\frac{n-1}{n(p+1)}}\leq \left({n(p+1)}{\delta}\right)^{\frac{1}{p+1}}\left(\int_{X}\varphi_{-}^{p+1}\frac{\omega_{\varphi}^{n}}{V}\right)^{\frac{1}{p+1}}.
\end{align*}
Now we start the standard iteration: Let $p_{0}:=1$ and $p_{j+1}+1:= \frac{n}{n-1}(p_{j}+1)$. Then we have that
\begin{align*}
\begin{split}
||\varphi_{-}||_{p_{j+1}+1}\leq C^{\frac{1}{p_{j}+1}}(p_{j}+1)^{\frac{1}{p_{j}+1}}||\varphi_{-}||_{p_{j}+1}&\leq \dots\\
&\leq C^{\sum_{i=0}^{\frac{1}{p_{i}+1}}}\prod_{i=0}^{j}(p_{i}+1)^{\frac{1}{p_{j}+1}}||\varphi_{-}||_{2}.
\end{split}
\end{align*}
That is to say
\begin{align*}
||\varphi_{-}||_{p_{j+1}+1}\leq C^{\sum_{i=0}^{\frac{1}{p_{i}+1}}}\prod_{i=0}^{j}(p_{i}+1)^{\frac{1}{p_{j}+1}}||\varphi_{-}||_{2}.
\end{align*}
It is not hard to check that the infinite product converges. Taking limits as $j\longrightarrow \infty$ gives
\begin{align*}
||\varphi_{-}||_{\infty}\leq C\left(\int_{X}\varphi_{-}^{2}\frac{\omega_{\varphi}^{n}}{V}\right)^{\frac{1}{2}}\leq
||\varphi_{-}||_{\infty}^{\frac{1}{2}}\ C\left(\int_{X}\varphi_{-}\frac{\omega_{\varphi}^{n}}{V}\right)^{\frac{1}{2}}\ .
\end{align*}
Which implies
\begin{align*}
||\varphi_{-}||_{\infty}\leq
 C^{2}\left(\int_{X}\varphi_{-}\frac{\omega_{\varphi}^{n}}{V}\right)\ .
\end{align*}
Since $\varphi_{t}\leq C$ as $t\rightarrow 0$ we have
\begin{align*}
-\mbox{inf}_{X}\varphi_{t}\leq C_{1} \int_{X}(-\varphi_{t})\frac{\omega_{\varphi}^{n}}{V}+C_{2} \ .
\end{align*}
Now, by the Green identity we deduce
\begin{align*}
\mbox{Osc}_{X}(\varphi_{t}):= \mbox{Sup}_{X}(\varphi_{t})-\mbox{Inf}_{X}(\varphi_{t})\leq C _{1}\left(\int_{X}\varphi_{t}\omega^{n}-
\int_{X} \varphi_{t}\omega_{\varphi}^{n} \right)+ C_{2}\ .
\end{align*}
Using the properness assumption gives:
\begin{align*}
f(\mbox{Osc}_{X}(\varphi_{t})) \leq \nu_{\omega}(\varphi_{t}).
\end{align*}
Now we are prepared to complete the proof of the corollary. \\
\noindent \textbf{Case 1}:\\
Assume that $X^{\lambda(0)}\neq X$ and moreover that $X^{\lambda(0)}$  is reduced, then by the same argument as in \cite{psc} we have
\begin{align*}
\mbox{lim}_{t \rightarrow 0}\mbox{Osc}_{X}(\varphi_{t})\rightarrow \infty.
\end{align*}
Consequently we deduce that
\begin{align*}
\mbox{lim}_{t \rightarrow 0}\nu_{\omega}(\varphi_{t})\rightarrow \infty .
\end{align*}
 corollary 1 yields the precise asymptotics \footnote{Recall that when $X^{\lambda(0)}$ is multiplicity free $\Psi(\lambda(t))=O(1)$.}
\begin{align*}
 \nu_{\omega}(\varphi_{\lambda(t)}) = F_{1}(\lambda)\log(t^2)+O(1) .
\end{align*}
This forces the desired sign  $F_{1}(\lambda)<0$.\\
 \textbf{Case 2}:\\
 If $X^{\lambda(0)}$ is nonreduced, then $\Psi(\lambda(t))\rightarrow -\infty$, however, under the properness assumption the K-Energy is bounded from below, and we again have that $F_{1}(\lambda)<0$. This completes the proof of corollary 2. $\Box$
 \bibliography{ref}
\end{document}